\documentclass{birkart}

\usepackage[dvips]{graphicx}

 \newtheorem{thm}{Theorem}[section]
 \newtheorem{cor}[thm]{Corollary}
 \newtheorem{lem}[thm]{Lemma}
 \newtheorem{prop}[thm]{Proposition}
 \theoremstyle{definition}
 
 \theoremstyle{remark}
 \newtheorem{rem}[thm]{Remark}
 
 \numberwithin{equation}{section}

\def\Xint#1{\mathchoice
   {\XXint\displaystyle\textstyle{#1}}%
   {\XXint\textstyle\scriptstyle{#1}}%
   {\XXint\scriptstyle\scriptscriptstyle{#1}}%
   {\XXint\scriptscriptstyle\scriptscriptstyle{#1}}%
   \!\int}
\def\XXint#1#2#3{{\setbox0=\hbox{$#1{#2#3}{\int}$}
     \vcenter{\hbox{$#2#3$}}\kern-.5\wd0}}

\def\dashint{\Xint-}
\def\fint{\dashint}

\def\rrho{\widetilde{\rho}}

\def\M{\mathcal{M}}
\def\Com#1{\mathbb{C}^{#1}}
\def\R#1{\mathbb{R}^{#1}}
\def\Z#1{\mathbb{Z}^{#1}}

\def\DD{\mathbb{D}}

\def\z{\zeta}

 \DeclareMathOperator{\re}{Re}

\DeclareMathOperator{\sgn}{sgn} 
\DeclareMathOperator{\supp}{supp}

\begin{document}
%
\title[Subharmonicity of higher dimensional exponential transforms]
 {Subharmonicity of higher dimensional\\ exponential transforms}
\author[Tkachev~ V.~G.]{Tkachev Vladimir G.}

\address{%
Volgograd State University, Department of Mathematics, 2-ya
Prodolnaya 30, 400062, Volgograd\\
Russia}

\email{vladimir.tkatchev@volsu.ru}

\thanks{The author was supported by grant RFBR
no.~03-01-00304 and by G\"oran Gustafsson Foundation.}

\subjclass{Primary 31C05; Secondary 44A12, 53C65}

\keywords{Riesz potential, exponential transform, subharmonic
function}

\date{\today{}}
\dedicatory{To Harold Shapiro on his 75th Anniversary, with
admiration.}

\begin{abstract}
Our main result states that the function $(1-E_\rho)^{(n-2)/n}$ is
subharmonic, where $0\leq \rho\leq 1$ is a density function in
$\R{n}$, $n\geq 3$, and $ E_\rho(x)=\exp\left(-\frac{2}{n}\;\fint
\frac{\rho(\zeta)d\zeta}{|\zeta-x|^n}\right)$,  is the exponential
transform of $\rho$. This answers in affirmative the recent
question posed by B. Gustafsson and M. Putinar in \cite{GP-prep}.

\end{abstract}

\maketitle


\section{Introduction}

The exponential transform can be viewed as a potential depending
on a domain in $\R{n}$, or more generally on a measure having a
\textit{density} function $\rho(x)$ (with compact support) in the
range $0\leq \rho\leq 1$. The two-dimensional version
\begin{equation}\label{intro}
E_\rho(z,w)=\exp\left[-\frac{1}{\pi}\int
\frac{\rho(\z)\;dA(\z)}{(\z-z)(\bar{\z}-\bar{w})}\right]
\end{equation}
has appeared in operator theory, as a determinantal-characteristic
function of certain close to normal operators \cite{Carey},
\cite{Pinkus}, and has previously been studied and proved to be
useful within operator theory, moment problems and other problems
of domain identification, and for proving regularity of free
boundaries (see \cite{GP-prep}, \cite{Put} for further
references). A corresponding exponential transform on the real
axis was already known and used by A.A.~ Markov (in the 19th
century) and later by N.I.~ Akhiezer and M.G.~ Krein in their
studies of one-dimensional moment problems \cite{Akh},
\cite{Akh-Kr}, (see, also \cite{Krein}).

In \cite{GP-prep} the diagonal version of (\ref{intro})
\begin{equation*}\label{intro1}
E_\rho(x)=\exp\left[-\frac{2}{\omega_n}\int
\frac{\rho(\zeta)d\zeta}{|x-\zeta|^{n}}\right]
\end{equation*}
is studied in higher dimensional case $n\geq 3$. Here $\omega_n$
denotes the $(n-1)$-dimensional Lebesgue measure of the unit
sphere in $\R{n}$.

Clearly, $0<E_\rho(x)<1$ for all $x\not\in\supp \rho$. In
particular, it was shown in \cite{GP-prep} that $E_\rho$ is a
subharmonic function. In two dimensions it is also known that
function $\ln (1 - E_\rho)$ is subharmonic, which is a stronger
statement. Here we extend the mentioned sub/superharmonicity in
dimension $n \geq 3$ thereby answering in affirmative a recent
question \cite[p.~566]{GP-prep}:

\begin{thm}\label{theo:harm1} Let  $E_\rho(x)$ be the
exponential transform of a density $\rho\not\equiv 0$. Then the
function
\begin{equation}\label{e-tilda1}
\left\{%
\begin{array}{ll}
    \ln (1-E_\rho), & \hbox{if \quad $n=2$,} \\
    \\
    \frac{1}{n-2}(1-E_\rho)^{(n-2)/n}, & \hbox{if \quad $n\geq 3$,} \\
\end{array}%
\right.
\end{equation}
is subharmonic outside $\supp \rho$.
\end{thm}

In fact, we show that a stronger version holds. To formulate it we
need some notation. Given an integer $n\geq 1$, we define $
\M_n(t)$ as the solution of the following ODE:
\begin{equation}\label{equ:diff}
\M_n'(t)=1-\M_n^{2/n}(t), \qquad \M(0)=0.
\end{equation}
We call $\M_n(t)$ the \textit{profile} function.

\begin{thm}\label{lem:harm-M}
For $n\geq 2$ let $\rho$ be a density function and
\begin{equation}\label{V}
V_\rho(x)=-\frac{n}{2}\ln E_\rho(x)\equiv \frac{n}{\omega_n}\int
\frac{\rho(\zeta)d\zeta}{|x-\zeta|^n}.
\end{equation}
Then the function
\begin{equation}\label{V-1}
\left\{
\begin{array}{ll}
    \log \M_2(V_\rho(x)), &  \text{if}\quad \hbox{$n=2$} \\
    \\
    \left[\M_n(V_\rho(x))\right]^{(n-2)/n}, &\text{if}\quad \hbox{$n\ne 2$} \\
\end{array}
\right.
\end{equation}
is subharmonic outside the support of $\rho$. Moreover, this
function is harmonic in $\R{n}\setminus B$, if $B$ is an arbitrary
Euclidean ball and $\rho=\chi_B$ is its characteristic function.
\end{thm}

We discuss properties of the profile function in more detail in
Section~\ref{sec:Mn}. In particular we show that $1-\M_n(x)$ is a
completely monotonic function in $\R{+}$.

\section{The main inequality}

\subsection{Variational problem}
Let $x=(x_1,y)\in\R{n}$, $y=(x_2,\ldots,x_n)$, and
$$
\R{n}_\pm=\{x=(x_1,y):\pm \,x_1>0\}.
$$
Given a measurable function $h(x)$ we denote by $\mathcal{J}(h)$
the integral
$$
\mathcal{J}(h)=\fint_{\R{n}}h(x)\;dx=\frac{n}{\omega_n}\int_{\R{n}}h(x)\;dx
$$
where $dx=dx_1dy$ denotes the $n$-dimensional Lebesgue measure in
$\R{n}$. In what follows we fix the following notations
$$
f(x)=\frac{1}{|x|^{n-2}}, \qquad g(x)=\frac{1}{|x|^{n}}, \qquad
\varphi(x)=\frac{x_1}{|x|^{n}},
$$
and suppose that $\rho(x)$ is a density function such that
$$
0\leq \rho(x)\leq 1.
$$
If $n=2$ we assume that $f(x)\equiv 1$. Throughout this section,
unless otherwise stated, we will assume that $\rho\ne0$ on a
non-null set and
 the support of $\rho$ does not contain a neighborhood of the
 origin. We write
\begin{equation}\label{equ:constrain}
\rho\in \mathcal{H}(w)\qquad\Leftrightarrow \qquad
\mathcal{J}(\rho g)\equiv
{\displaystyle\fint_{\R{n}}\frac{\rho}{|x|^{n}} \;dx}=w\geq 0.
\end{equation}

Our main subject is the ratio
$$
\Phi(\rho)=\frac{\mathcal{J}^2(\varphi\rho)}{\mathcal{J}(f\rho)}.
$$

\begin{thm}\label{theo:1}
Let $\rho$ be a density function, $0\not\in\supp\rho$. Then
\begin{equation}\label{e:equal}
\max_{\rho\in \mathcal{H}(w)}\Phi(\rho)=\M_n(w).
\end{equation}
For any $w>0$ the maximum is attained when $\rho(x)$ is the
characteristic function of the ball centered at
$(\tau,\mathbf{0})$ of radius $\tau \M_n(w)^{1/n}$, with $\tau>0$.
\end{thm}


We mention two limit cases of the last assertion. Namely, the
boundedness of maximum in (\ref{e:equal}) easily follows from
$\varphi^2\leq fg$ and the Cauchy-Schwarz inequality:
\begin{equation}\label{equ:Cauchy}
\frac{\mathcal{J}^2(\varphi\rho)}{\mathcal{J}(f\rho)}\leq
\mathcal{J}(g\rho)=w.
\end{equation}
On the other hand, it was shown by Gustafsson and Putinar in
\cite[p.~563]{GP-prep} that
\begin{equation}\label{equ:GP1}
\frac{\mathcal{J}^2(\varphi\rho)}{\mathcal{J}(f\rho)}< 1
\end{equation}
does hold. The last means that  inequality (\ref{equ:GP1})
considerably refines (\ref{equ:Cauchy}) when $w>1$ while the first
estimate becomes to be sharper when $w$ is a small value.

\begin{cor}For any density function $\rho(x)$,
$0\not\in\supp\rho$, the following sharp inequality holds
\begin{equation}\label{var}
\left(\fint_{\R{n}}\frac{x_1\rho(x)}{|x|^n}\; dx\right)^2\leq
\M_n\left(\fint_{\R{n}}\frac{\rho(x)}{|x|^n}\; dx\right)
{\displaystyle\fint_{\R{n}}\frac{\rho(x)}{|x|^{n-2}} \;dx}
\end{equation}
\end{cor}

The inversion $x\to x/|x|^2$ gives another equivalent form of the
preceding property

\begin{cor}\label{cor:invers}For any density function $\rho(x)$,
$0\not\in\supp\rho$, the following sharp inequality holds
\begin{equation}\label{var+1}
\left(\fint_{\R{n}}\frac{x_1\rho(x)}{|x|^{n+2}}\; dx\right)^2\leq
\M_n\left(\fint_{\R{n}}\frac{\rho(x)}{|x|^n}\; dx\right)
{\displaystyle\fint_{\R{n}}\frac{\rho(x)}{|x|^{n+2}} \;dx}
\end{equation}
\end{cor}

\begin{rem}
We note that for $n\geq 3$ the above inequality (\ref{var}) can be
interpreted as a pointwise estimate on the  Coulomb potential
$$
U_\rho(x)=\fint \frac{\rho(\zeta) dx}{|x-\zeta|^{n-2}}
$$
with an \textit{bounded} density function $\rho$, $0\leq \rho\leq
1$. Indeed, using the inversion in $\R{n}$ we see that (\ref{var})
is equivalent to
\begin{equation*}\label{equ:potential}
|\nabla U_\rho(x)|^2\leq \M_n[V_\rho(x)] U_\rho(x),\qquad x
\not\in \supp \rho,
\end{equation*}
where $V_\rho(x)$ is defined by (\ref{V}). In particularly,
$\M_n(w)<1$ gives us the inequality due to Gustafsson and Putinar
\cite{GP-prep}:
\begin{equation*}\label{equ:GP}
|\nabla U_\rho(x)|^2< U_\rho(x), \qquad x \not\in \supp \rho.
\end{equation*}
\end{rem}

\subsection{Auxiliary integrals}\label{sec:int}
In order to prove Theorem~\ref{theo:1}, we need to evaluate the
integrals in (\ref{var}) for a specific choice of the density
function. Namely, let $\tau>\alpha>0$ and consider the following
density function
$$
\widehat{\rho}(x)=\chi_{\mathbb{D}}(x),
$$
where
\begin{equation}\label{D-a}
\mathbb{D}\equiv\mathbb{D}(\alpha,\tau):=\biggl\{x=(x_1,y):\quad
(x_1-\tau)^2+|y|^2<\tau^2-\alpha^2\biggr\}.
\end{equation}

First, we note that the function $f(x)=|x|^{2-n}$ is harmonic in
$\overline{\mathbb{D}}$. Using the fact that the ball $\mathbb{D}$
is of radius $\sqrt{\tau^2-\alpha^2}$ and centered at
$x=(\tau,\mathbf{0})$, we have by the mean value theorem
\begin{equation}\label{equ:f}
\mathcal{J}(\rrho f)=\fint_{\mathbb{D}}\frac{dx}{|x|^{n-2}}
=\frac{({\tau^2-\alpha^2})^{n/2}}{\tau^{n-2}}
=\alpha^2\frac{\sinh^n\xi}{\cosh^{n-2}\xi}
\end{equation}
where
\begin{equation}\label{equ:xi}
\cosh\xi=\frac{\tau}{\alpha}.
\end{equation}

Similarly, harmonicity of $\varphi(x)=x_1|x|^{-n}$ implies
\begin{equation*}\label{equ:phi}
\mathcal{J}(\rrho \varphi)=\alpha
\frac{\sinh^n\xi}{\cosh^{n-1}\xi}.
\end{equation*}

To evaluate $\mathcal{J}(\rrho g)$ we consider the following
auxiliary  function
$$
\lambda(x)=\frac{|x|^2+\alpha^2}{2\tau x_1}.
$$
Then $\lambda(x)$ is positive on $\mathbb{D}$ and ranges in
$$
\frac{\alpha}{\tau}\leq \lambda(x)<1, \quad x\in \mathbb{D}.
$$
Moreover, it is easy to see that
\begin{equation}\label{e:d-e}
\qquad \lambda(x)\equiv z, \qquad x\in
S(z)=\partial\mathbb{D}(\alpha,\tau z).
\end{equation}
Hence, the co-area formula yields
\begin{equation}\label{equ:co}
\mathcal{J}(\rrho g)=\fint_{\mathbb{D}}\frac{dx}{|x|^n}=
\frac{n}{\omega_n}\int_{\alpha/\tau}^1dz
\int_{S(z)}\frac{dS}{|x|^n|\nabla \lambda(x)|}.
\end{equation}
Here $dS$ is the $(n-1)$-dimensional surface measure of the level
set $S(z)$.

On the other hand, we have for the gradient
$$
|\nabla \lambda|^2=\frac{|y|^2}{\tau^2 x_1^2}+
\frac{(x_1^2-\alpha^2-|y|^2)^2}{4\tau^2 x_1^4}
$$
which by virtue of (\ref{e:d-e}) implies the corresponding value
on the level set $S(z)$:
$$
|\nabla \lambda|^2\biggr|_{S(z)}=\frac{\tau^2z^2-\alpha^2}{\tau^2
x_1^2}.
$$
Substitution of the last expression into (\ref{equ:co}) yields
$$
\mathcal{J}(\rrho
g)=\frac{n}{\omega_n}\int_{\alpha/\tau}^1\frac{\tau
dz}{\sqrt{\tau^2z^2-\alpha^2}} \int_{E(z)}\frac{x_1}{|x|^n}dS.
$$
Since  $\varphi(x)=x_1|x|^{-n}$ in the inner integral is a
harmonic function and $S(z)$ is a sphere, we have by the mean
value theorem
\begin{equation*}
\begin{split}
\mathcal{J}(\rrho
g)&=\frac{n}{\omega_n}\int_{\alpha/\tau}^1\frac{\tau
dz}{\sqrt{\tau^2z^2-\alpha^2}}\cdot\frac{(\tau^2z^2-\alpha^2)^{(n-1)/2}}{(\tau
z)^{n-1}}=n\int_{0}^{\xi}\tanh^{n-1} t dt.\\
\end{split}
\end{equation*}
where $\xi$ is defined by (\ref{equ:xi}). Thus we obtain
\begin{equation}\label{equ:co1}
\mathcal{J}(\rrho g)=T_n(\tau/\alpha)=T_n(\xi):=
n\int_{0}^{\xi}\tanh^{n-1}t\; dt.
\end{equation}
We point out that \textit{the latter integral depends only on the
ratio} $\tau/\alpha$. One can easy verify that
\begin{equation}\label{new-1}
\M_n(T_n(\xi))\equiv
\tanh^n\xi=\left(\frac{\sqrt{\tau^2-\alpha^2}}{\tau}\right)^n.
\end{equation}

\begin{rem}\label{rem:M}
After a suitable shift in the $x_1$-direction, the last
computation is equivalent to the following relation
\begin{equation}\label{equ:co1y}
\M_n\left(\fint_{\mathbb{B}(R)}\frac{d\zeta}{|x-\zeta|^{n}}\right)
=\left(\frac{R}{|x|}\right)^n,
\end{equation}
which holds for any ball  $\mathbb{B}(R)$  of radius $R$ centered
at the origin.
\end{rem}

\subsection{Proof of Theorem~\ref{theo:1}}
Let us denote
\begin{equation}\label{equ:ex+}
M_n^+(w)=\sup_{\mathrm{\rho\in\mathcal{R}}}\Phi(\rho)
\end{equation}
where $\mathcal{R}$ denotes  the class of all density functions
$\rho$ such that $\supp \rho\cap \R{n}_-$ has null measure. Then
Theorem~\ref{theo:1} follows from the following lemmas.

\begin{lem}\label{lem:1}
$M_n^+(w)=\M_n(w)$.
\end{lem}

\begin{lem}\label{lem:2}
$ \sup_{\mathrm{\rho}}\Phi(\rho)=M_n^+(w). $
\end{lem}

\begin{proof}[Proof of Lemma~\ref{lem:1}]
Our first step is to reduce the problem (\ref{equ:ex+}) to the
following linear extremal problem with additional constraints:
\begin{equation}\label{equ:extt}
N_n(w):=\sup_{\rho\in\mathcal{R}}\{\mathcal{J}(\rho \varphi):
\quad \mathcal{J}(\rho f)=1, \; \mathcal{J}(\rho g)=w\}.
\end{equation}
Then we have
\begin{equation}\label{equ:NM}
M_n^+(w)=N_n^2(w).
\end{equation}
Indeed, in order to prove (\ref{equ:NM}), let $\rho_a(x)=\rho(ax)$
be a homothety of $\rho(x)$ with positive coefficient $a$.
Clearly, this transformation preserves the class $\mathcal{R}$. On
the other hand, one can easily see that
$$
\Phi(\rho_a)=\Phi(\rho)
$$
by the virtue of homogeneity of $\Phi$. Moreover,
$$
\mathcal{J}(\rho_a \varphi)=\frac{1}{a}\mathcal{J}(\rho \varphi),
\qquad \mathcal{J}(\rho_a f)=\frac{1}{a^2}\mathcal{J}(\rho f),
$$
which proves (\ref{equ:NM}).

Next, we claim that for any nonnegative $w$ there exists an
$\alpha>0$ and $\tau>\alpha$ such that
\begin{equation}\label{equ:main}
\mathcal{J}(\rrho  f)=1, \qquad \mathcal{J}(\rrho  g)=w,
\end{equation}
where $\rrho=\chi_{\mathbb{D}(\alpha,\tau)} $ is the
characteristic function of the ball $\mathbb{D}(\alpha,\tau)$ in
(\ref{D-a}). Indeed, using the definition of function $T_n(t)$ in
(\ref{equ:co1}) one can easily see that  there exist a unique root
$\xi>0$ of the following equation
\begin{equation}\label{Tw}
T_n(\xi)=w.
\end{equation}
Then we chose $\alpha>0$ such that
$$
\alpha^2=\frac{\cosh^{n-2}\xi}{\sinh ^{n}\xi},
$$
and let $\tau =\alpha\cosh\xi$. Now (\ref{equ:main}) immediately
follows from  (\ref{equ:f}) and (\ref{equ:co1}).

Thus, the function $\rrho (x)$ satisfies (\ref{equ:main}) and it
follows that it is admissible for the problem (\ref{equ:extt}).
This implies
$$
N_n(w)\geq \mathcal{J}(\rrho \varphi).
$$

To prove that the inverse inequality holds, we fix any function
$\rho\in\mathcal{R}$ which is admissible for (\ref{equ:extt}).
Then
$$
\mathcal{J}(\rrho(f+\alpha^2 g))=\mathcal{J}(\rho(f+\alpha^2
g))=1+\alpha^2w.
$$
The last property means that both the functions $\rho$ and
$\rrho$\; are test functions for the following extremal problem
\begin{equation}\label{e:addd}
\sup_{\rho\in\mathcal{R}}\{\mathcal{J}(\rho\varphi):\quad
\mathcal{J}(\rho (f+\alpha^2 g))=1+w\alpha^2\}.
\end{equation}
Let us consider the ratio
$$
h(x):=\frac{\varphi(x)}{f(x)+\alpha^2 g(x)}=
\frac{x_1}{|x|^2+\alpha^2}.
$$
Then,
$$
{\{x\in\R{n}:h(x)>\frac{1}{2\tau}\}}=\mathbb{D}(\alpha,\tau),
$$
and it  follows from the Bathtub Principle \cite[p.~28]{Lieb} that
$\rrho$ is the extremal density for (\ref{e:addd}). Thus, we have
$$
\mathcal{J}(\rho\varphi)\leq \mathcal{J}(\rrho\varphi),
$$
and consequently
$$
N_n(w)\leq \mathcal{J}(\rrho\varphi).
$$
Hence, we conclude that
$$
N_n(w)= \mathcal{J}(\rrho\varphi)=\alpha
\frac{\sinh^n\xi}{\cosh^{n-1}\xi}.
$$

Now, it follows from (\ref{equ:NM}) and our choice of $\alpha$
that
$$
M^+_n(w)=N_n^2(w)= \alpha^2
\frac{\sinh^{2n}\xi}{\cosh^{2n-2}\xi}=\tanh^n\xi,
$$
and from (\ref{new-1}), we find
$$
M^+_n(w)=\M_n(T_n(\xi))=\M_n(w),
$$ and the lemma follows.
\end{proof}

\begin{proof}[Proof of Lemma~\ref{lem:2}]
It suffices only to prove the one-side inequality
\begin{equation}\label{equ:+}
\sup_{\mathrm{\rho}}\Phi(\rho)\leq M_n^+(w).
\end{equation}

Let $\rho$ is an arbitrary admissible for (\ref{equ:constrain})
density function. Excluding the trivial case $\rho\in\mathcal{R}$
we distinguish two rest cases:

(i)  the set $\supp  \;\rho\cap \R{n}_+$ has the null measure;

(ii) the set $\supp  \;\rho$ has non-zero counterpart in the both
half-spaces.

Let $\rho$ satisfies (i). Then the function
$$
\rho^*(x_1,y):=\rho(-x_1,y)
$$
belongs to $\mathcal{R}$, and it follows that
\begin{equation}\label{M+}
\sup_{\mathrm{\rho\in
\mathrm{(i)}}}\Phi(\rho)=\sup_{\mathrm{\rho\in\mathcal{R}}}\Phi(\rho)
=M_n^+(w).
\end{equation}

Now, let $\rho$ satisfies (ii). We set $
\rho^{\pm}(x)=\chi_{\R{n}_{\pm}}(x)\rho(x)$. Then
$$
\mathcal{J}(\rho \varphi)=\mathcal{J}(\rho^+
\varphi)-\mathcal{J}((\rho^-)^* \varphi),
$$
$$
\mathcal{J}(\rho f)=\mathcal{J}(\rho^+ f)+\mathcal{J}((\rho^-)^*
f),
$$
where the last integrals are positive. Using an elementary
inequality
$$
\frac{(a-b)^2}{c+d}\leq \max
\left[\frac{a^2}{c},\frac{b^2}{d}\right]
$$
which holds for any set of positive numbers $a,b,c,d$, we conclude
that
$$
\Phi(\rho)=\frac{\mathcal{J}^2(\rho \varphi)}{\mathcal{J}(\rho
f)}\leq \max [\Phi(\rho^+),\Phi((\rho^-)^*) ].
$$
Hence, we have by Lemma~\ref{lem:1}
$$
\Phi(\rho)\leq \max [M^+_n(w_1), M^+_n(w_2)]=\max [\M_n(w_1),
\M_n(w_2)],
$$
where
$$
w_1=\mathcal{J}(\rho^+ g), \qquad w_2=\mathcal{J}((\rho^+)^* g).
$$
But
$$
w=\mathcal{J}(\rho g)=w_1+w_2,
$$
whence $w_i\leq w$, $i=1,2$. Since $\M_n$ is an increasing
function we obtain $\Phi(\rho)\leq \M_n(w)$, and consequently
$$
\sup_{\mathrm{\rho\in
\mathrm{(ii)}}}\Phi(\rho)\leq\M_n(w)=M_n^+(w).
$$

Combining the last inequality with (\ref{M+}) we obtain
$$
\sup_{\mathrm{\rho}}\Phi(\rho)=\sup_{\mathcal{R}\cup\mathrm{(i)\cup
\mathrm{(ii)}}}\Phi(\rho) \leq M_n^+(w)
$$
which proves (\ref{equ:+}).
\end{proof}

\section{Proof of the main results}

\begin{lem}\label{lem:exp}
For any $n\geq 1$ we have
\begin{equation}\label{equ:super}
\M_n(w)\leq Q_n(w):= \frac{e^{2w/n}-1}{e^{2w/n}-\frac{n-2}{n}}.
\end{equation}
\end{lem}

%

\begin{proof}
Note that in the cases $n=1,2$, we  have
$$
\M_1(w)=\tanh w= \frac{e^{2w}-1}{e^{2w}+1},
$$
$$
\M_2(w)=1-e^{-w}
$$
which turns (\ref{equ:super}) into equality.

Now, let $n\geq 3$. We have $M_n(0)=Q_n(0)=0$ and by the
definition (\ref{equ:diff}) it suffices only to prove that
\begin{equation}\label{equ:Q}
Q'_n(w)\geq 1-Q_n^{2/n}, \quad w>0.
\end{equation}

We have
$$
Q_n'(w)=(1-\frac{n-2}{n}Q_n)(1-Q_n)
$$
and (\ref{equ:Q}) becomes to be equivalent to the inequality
$$
\frac{1-t^{1-\gamma}}{1-t}< 1-\gamma t,
$$
where $t=Q_n(w)\in (0,1)$ and $\gamma=(n-2)/n$. To verify the last
inequality we rewrite it in the form
$$
\frac{1-t^\gamma}{1-t}> \gamma t^\gamma.
$$
For $t\in (0,1)$, the function in the left hand side is a
decreasing function  while the right hand side member is an
increasing one. Since the both functions have the same limit value
$\gamma$ at $t=1$, we have the desired inequality.
\end{proof}

\begin{proof}[Proof of Theorem~\ref{theo:harm1}]
Let $f(x)$ denote the function in (\ref{e-tilda1}). Then we have
for any $n\geq 2$ and $x\not\in\supp\rho$
\begin{equation*}
\begin{split}
\nabla f(x)&= -(1-E_\rho)^{-2/n}\nabla E_\rho,
\\
\Delta
f(x)&=-\frac{2}{n}(1-E_\rho)^{-\frac{2+n}{n}}\biggl[\frac{n}{2}(1-E_\rho)\Delta
E_\rho+|\nabla E_\rho|^2\biggr].
\end{split}
\end{equation*}
Then the inequality $\Delta f(x)\geq 0$ to be proved becomes
\begin{equation}\label{becomes}
\frac{n}{2}(1-E_\rho)\Delta E_\rho+|\nabla E_\rho|^2\leq 0.
\end{equation}
On the other hand,
\begin{equation*}
\begin{split}
\nabla E_\rho(x)&=2E_\rho(x) \fint
\frac{(x-\zeta)\rho(\z)d\z}{|x-\zeta|^{n+2}} ,
\\
\Delta E_\rho(x)&=4E_\rho(x)\biggl( \biggl|\fint
\frac{(x-\zeta)\rho(\z)d\z}{|x-\zeta|^{n+2}}\biggr|^2-\fint
\frac{\rho(\z)d\z}{|x-\zeta|^{n+2}} \biggr),
\end{split}
\end{equation*}
and (\ref{becomes}) becomes
\begin{equation}\label{becomes2}
\left(1-\frac{n-2}{n}E_\rho\right)\biggl|\fint
\frac{(x-\zeta)\rho(\z)d\z}{|x-\zeta|^{n+2}}\biggr|^2\leq
(1-E_\rho)\fint \frac{\rho(\z)d\z}{|x-\zeta|^{n+2}}.
\end{equation}
In order to prove (\ref{becomes2}) we can assume without loss of
generality that $x=0$. In this case, after a suitable rotation we
can write the vector integral as follows
$$
\biggl|\fint \frac{\zeta\rho(\z)d\z}{|\zeta|^{n+2}}\biggr|= \fint
\frac{\zeta_1\rho(\z)d\z}{|\zeta|^{n+2}}.
$$
Thus, we arrive at the inequality to be proved
\begin{equation}\label{becomes3}
\left(1-\frac{n-2}{n}e^{-\frac{2w}{n}}\right)\biggl(\fint
\frac{\zeta_1\rho(\z)d\z}{|\zeta|^{n+2}}\biggr)^2\leq
(1-e^{-\frac{2w}{n}})\fint \frac{\rho(\z)d\z}{|\zeta|^{n+2}},
\end{equation}
with
$$
w=\fint \frac{\rho(\z)d\z}{|\zeta|^{n}}.
$$
But, it is easy to see that (\ref{becomes3}) follows from
Corollary~\ref{cor:invers} and Lemma~\ref{lem:exp}. The theorem
follows.
\end{proof}

\begin{proof}[Proof of Theorem~\ref{lem:harm-M}]
Let $F(x)$ denote the function in (\ref{V-1}) and
$V(x)=V_\rho(x)$. Then the argument similar to that above yields
for $n\geq 3$
\begin{equation}\label{equ:ii}
\begin{split}
\Delta F(x) &=\Delta (\M_n(V))^{(n-2)/n}\\
&=(1-\M^{2/n}_n(V))\left(\frac{n-2}{n}
\M^{-\frac{2}{n}}_n(V)\Delta V -
\frac{2(n-2)}{n^2}\M_n^{-\frac{2+n}{n}}(V) |\nabla V|^2\right)\\
&=2(n-2)(1-\M^{2/n}_n(V))\biggl[\M_n(V)B - |A|^2\biggr],
\end{split}
\end{equation}
where
$$
A=\fint \frac{(x-\zeta)\rho(\zeta)d\zeta}{|x-\zeta|^{n+2}}, \qquad
B=\fint \frac{\rho(\zeta)d\zeta}{|x-\zeta|^{n+2}}.
$$
Similarly, we have for $n=2$
$$
\Delta F(x)=\frac{1-\M_2(V)}{\M^2_2(V)}\biggl[\M_2(V)B -
|A|^2\biggr].
$$
Hence, for all integer $n\geq 2$, the sign of the Laplacian
$\Delta F(x)$ coincides with the sign of $[M_n(V)B - |A|^2]$.

Let us fix an arbitrary point $x\not\in\supp\rho$. Then after a
suitable rotation we can reduce the vector integral $A$ to the
scalar one such that the value in last brackets in (\ref{equ:ii})
becomes
$$
\M_n\left(\fint \frac{\rho_1(\zeta)d\zeta}{|\zeta|^{n}}\right)
\fint \frac{\rho_1(\zeta)d\zeta}{|\zeta|^{n+2}}- \left(\fint
\frac{\zeta_1\rho_1(\zeta)d\zeta}{|\zeta|^{n+2}}\right)^2,
$$
where $\rho_1(\zeta)$ is the correspondent transformed density.
Then Corollary~\ref{cor:invers} again implies that the latter
difference is nonnegative and subharmonicity of $\mathcal{E}_\rho$
easily follows.

Now, let us prove the second assertion of the theorem.  Let
$\mathbb{B}(R)$ be the  ball of radius $R$ with center at the
origin and $\widehat{\rho}(x)=\chi_{\mathbb{B}(R)}(x)$ be the
corresponding characteristic function. Then
$$
V_{\widehat{\rho}}(x):=\fint_{\mathbb{B}(R)}
\frac{d\zeta}{|x-\zeta|^n},
$$
and  we have from (\ref{equ:co1y}) that in this case
$$
\M_n(V_{\widehat{\rho}}(x))=\left(\frac{R}{|x|}\right)^n,
$$
which obviously yields harmonicity of
$$
\left[\M_n(V_{\widehat{\rho}}(x))\right]^{(n-2)/n}=R^{n-2}|x|^{2-n}
$$
for $n\geq 3$, and
$$
\ln \M_2(V_{\widehat{\rho}}(x))=2\ln \frac{R}{|x|},
$$
if $n=2$. The theorem is completely proved.
\end{proof}

\section{The profile function}
\label{sec:Mn} Here we study the profile function $\M_n$ in more
detail. This higher transcendental function, apart of its
appearance in the above theorems, admits also number-theoretical
applications (e.g., in connection with the Euler-Mascheroni
constant $\gamma$, see Section~\ref{eul}). Our main result
(Theorem~\ref{theo:com} below) states that $1-M_n(w)$ is a
completely monotonic function. We also show
(Theorem~\ref{theo:phi}) that this function can be analytically
extended across $w=+\infty$
 by making use of a specific logarithmic
transformation.

\subsection{Complete monotonicity}
It is convenient to consider the general case of (\ref{equ:diff}).
Namely, given a real $\alpha>0$ we define $F_\alpha(x)$ as a
solution to the following ODE
\begin{equation}\label{equ:extM}
F_\alpha'(x)=1-F_\alpha^\alpha(x), \qquad F_\alpha(0)=0.
\end{equation}
Then for an integer $n$ we have
$
\M_n(w)=F_{2/n}(w).
$

We recall that a function $f(x)$ defined on $[0;+\infty)$ is said
to be \textit{completely monotonic} if
$$
(-1)^kf^{(k)}(x)\geq 0, \qquad x\in \R{+}.
$$

\begin{thm}
\label{theo:com} Let $\alpha>0$. Then

(i) $F_\alpha(x)$ is an increasing function for $x\geq 0$ such
that $F_\alpha(x):\R{+}\to [0;1)$;

(ii) for all $\alpha\in(0;1]$ the function
\begin{equation*}\label{e:tilde}
\widetilde{F}_\alpha(x)=1-F_\alpha(x)
\end{equation*}
is completely monotone on $\R{+}$.
\end{thm}

It follows from the well-known Bernstein's theorem \cite{Bern-AMF}
(see also \cite[p.~161]{Widder}) that $\widetilde{F}_\alpha(x)$ is
a Laplace transform of a positive measure supported on $\R{+}$.

\begin{cor}
\label{cor:bern} For all $\alpha\in(0;1]$ the following
Laplace-Stieltjes representation holds
\begin{equation}\label{equ:rr}
\widetilde{F}_\alpha(x)=\int\limits_{0}^{+\infty}e^{-xt}\;
d\sigma_{\alpha}(t)
\end{equation}
where $d\sigma_\alpha$ is a positive probability measure with
finite variation
\begin{equation}\label{equ:sigma}
\int\limits_{0}^{+\infty}d\sigma_{\alpha}(t)=\widetilde{F}_\alpha(0)=1.
\end{equation}
\end{cor}


The following subadditive property is a consequence of the general
result due to Kimberling \cite{Kimb} and concerns complete
monotonic functions satisfying (\ref{equ:sigma}).

\begin{cor}
\label{cor:add} For all $0<\alpha\leq 1$ the function
$\widetilde{F}_\alpha(x)$ is subadditive in the sense that
\begin{equation}\label{equ:add}
\widetilde{F}_\alpha(x)\widetilde{F}_\alpha(y)\leq
\widetilde{F}_\alpha(x+y).
\end{equation}
\end{cor}

\begin{rem}
It is easy to verify that for $\alpha>1$ the third derivative of
$F'''_\alpha(x)$ has no constant sign on $\R{+}$. Thus our
constraint is optimal for positive values of $\alpha$. On the
other hand, if $\alpha=1$ then $F_1(x)$ can be derived as follows
$$
F_1(x) =1-e^{-x}, \qquad \widetilde{F}_1(x)=e^{-x}
$$
which implies the complete monotonicity immediately. Moreover, in
the latter case $\widetilde{F}_1(x)$ satisfies a full additive
property instead of (\ref{equ:add}). We notice also that in this
case one can easily find that $d\sigma_1(t)=\delta_1(t)$ the
delta-Dirac probability measure supported at $t=1$. More
precisely, we have
$$
\sigma_1(t)=\chi_{[1,+\infty)}(t).
$$
\end{rem}

\begin{proof}[Proof of Theorem~\ref{theo:com}]
The only non-trivial part of the theorem is (ii). We notice first
that
$$
\widetilde{F}_\alpha(x)\geq 0, \qquad
\widetilde{F}^{(k)}_\alpha(x)=-F^{(k)}_\alpha(x), \quad
k=1,2,\ldots
$$
and
$$
F_\alpha''(x)=-\alpha(1-F_\alpha^\alpha)\frac{F_\alpha^\alpha}{F_\alpha(x)}.
$$

On the other hand, one can easily show by induction that the
following property holds for all $k\geq 0$
\begin{equation}\label{equ:rep}
F^{(k+2)}_\alpha(x)=\alpha t(1-t)\frac{H_k(t)}{F_\alpha(x)^{k+1}}
\end{equation}
where
$$
t=F_\alpha^\alpha(x)
$$
and $H_{j}(t)$ is a polynomial of degree at most $j$. Moreover, we
have the following recurrent relationship
\begin{equation}\label{equ:H}
H_{k+1}(t)=[(k+1-2\alpha)t-(k+1-\alpha)]H_k(t)+\alpha
t(1-t)H'_k(t), \qquad k\geq 2
\end{equation}
with initial condition
\begin{equation}\label{equ:H0}
H_{0}(t)=-1.
\end{equation}

 Since $t=F_\alpha^\alpha(x)$ ranges in
$[0;1)$ we have only to prove that the polynomials
$(-1)^{k+1}H_k(t)$ are nonnegative in $\Delta=[0,1)$.

We will use the following Bernstein-type transformation
$$
P^*(z)=(1+z)^nP\left(\frac{1}{1+z}\right), \qquad n\geq\deg P
$$
which transforms a polynomial $P$ to a polynomial of degree at
most $n$.

Let
$$
P(t)=a_0+a_1t+\ldots +a_nt^n
$$
(here we use the assumption that $\deg P\leq n$ and some
coefficients may vanish). Then we can write
\begin{equation}\label{equ:Berr}
P(t)=\sum_{j=0}^{n}b_jt^{n-j}(1-t)^{j}
\end{equation}
where
$$
P^*(z)=b_0+b_1z+\ldots +b_nz^n, \qquad z=\frac{1-t}{t}.
$$
We recall that (\ref{equ:Berr}) is the Bernstein-type expansion of
$P$ by the basis $t^{j}(1-t)^{n-j}$.

It follows then from (\ref{equ:Berr}) that if all (non-zero)
coefficients of the associate polynomial $P^*(z)$ have the same
sign: $\sgn b_j=\varepsilon$, then $P(t)$ changes no sign in
$\Delta$ and its sign coincides with $\varepsilon$.

Let $H_k^\star(z)$ be the associative polynomial for $H_k(t)$.
Then
$$
H_k(t)=t^kH^*_k\left(\frac{1-t}{t}\right)
$$
and
$$
H'_k(t)=kt^{k-1}H^*_k\left(\frac{1-t}{t}\right)-
t^{k-2}{H^*_k}'\left(\frac{1-t}{t}\right).
$$
It follows from (\ref{equ:H}) that
\begin{equation}\label{equ:star}
-H^*_{k+1}(z)=[\alpha+(k+1)(1-\alpha)z]H^*_k(z)+\alpha
z(1+z){H^*_k}'(z).
\end{equation}

We notice that by (\ref{equ:H0})
$$
H^*_0=H_0=-1.
$$
On the other hand, since $0\leq \alpha\leq 1$ the multipliers
$(\alpha+(k+1)(1-\alpha)z)$ and $\alpha z(1+z)$ in
(\ref{equ:star}) have positive coefficients with respect to $z$.
Hence, it immediately follows from (\ref{equ:star}) by induction
that all coefficients of $-H^*_{k+1}(z)$ have the same sign as
$H^*_k(z)$ does. Moreover, the sign of the coefficients of
$H^*_k(z)$ is $(-1)^{k+1}$ which yields by the above remark  that
$$
(-1)^{k+1}H_k(t)\geq 0, \qquad t\in\Delta.
$$
Clearly, the last property together with (\ref{equ:rep}) yields
the desired assertion.
\end{proof}

\subsection{Exponential series for the profile function}\label{eul}
Here we establish an explicit form of the above exponential
representation for $\M_n(x)$. As above, it is convenient to
consider a general $F_\alpha(x)$ instead of $\M_n(x)$
 (see the definition (\ref{equ:extM})).

Let
$$
\phi_\alpha(t):=1-{F}_\alpha\left(-\frac{1}{\alpha}\ln t\right).
$$
According to its definition, $\phi_\alpha(t)$ is defined in
$(0,1]$. But it turns out that a stronger property holds

\begin{thm}\label{theo:phi}The following properties hold:
\begin{enumerate}
\item[(i)]For any $\alpha>0$ the function $\phi_\alpha(t)$ admits
an analytic continuation on $(-\epsilon,1)$ with some $\epsilon>0$
depending on $\alpha$.

\item[(ii)] The corresponding Taylor series at $t=0$ are
\begin{equation}\label{e:Taylor}
\phi_\alpha(t)=\sum_{k=1}^{\infty}\sigma_k (\gamma_\alpha t)^k,
\end{equation}
where
\begin{equation*}\label{e:Taylor-gamma}
\begin{split}
\gamma(\alpha)&= \frac{1}{\alpha}\exp\biggl(-\int\limits_{0}^1
\frac{1-x^{\frac{1-\alpha}{\alpha}}}{1-x\phantom{mn}}\;dx\biggr),\\
\end{split}
\end{equation*}
and $\sigma_k$ are the coefficients defined by the following
recurrence
\begin{equation}\label{e:Taylor-coef}
\sigma_1:=1, \quad
\sigma_k=\frac{1}{k(k-1)}\sum_{\nu=1}^{k-1}\sigma_\nu
\sigma_{k-\nu} [(1+\alpha)\nu-\alpha k]\nu.
\end{equation}
\item[(iii)] If $\alpha\in (0,1)$ then $\sigma_k>0$ for all $k\geq
1$ and series (\ref{e:Taylor}) converges in $(-1,1)$.

\item[(iv)] For all $0<\alpha<1$, $\phi_\alpha(t)$ is a strictly
increasing convex function in $(-\infty,1)$.
\end{enumerate}
\end{thm}

\begin{rem}
The exact value of $\gamma_\alpha$ has the following form
\begin{equation}\label{gamma}
\ln \gamma_{\alpha}=-\Psi \left( 1/\alpha \right) -\gamma+\ln
\left( 1/\alpha \right),
\end{equation}
where $\Psi(z)$ is the Digamma function:
$\Psi(z)=\Gamma'(z)/\Gamma(z)$, and $\gamma=0.5772156\ldots$ is
the Euler-Mascheroni constant. The assertion of the theorem is
still valid for $\alpha=0$ which formally corresponds to
$n=\infty$. In this case, $\phi_0(x)$ satisfies the following ODE:
\begin{equation*}\label{e:ln}
\phi_0'(x)=-\frac{\ln(1-\phi_0(x))}{x}, \qquad \phi_0(0)=0.
\end{equation*}
It follows from (\ref{gamma}) that in this case
$\gamma_0=e^{\gamma}$.
\end{rem}

\begin{cor}
Let $n\geq 2$ be an integer. Then
\begin{equation*}\label{e:exp-rep}
1-\M_n(x)=\sum_{k=1}^{\infty}a_k e^{-2k x/n},
\end{equation*}
where $a_k=\sigma_k \gamma_{2/n}^k>0$ and the series converges for
all $x\geq 0$. In particular, the measure in (\ref{equ:rr}) is an
atomic measure supported at the set $\frac{2}{n} \Z{+}$.
\end{cor}

We are grateful to Bj\"orn Gustafsson for pointing out another
useful consequence of the preceding property. Let us define an
($n$-dimensional) version of the exponential transform as follows
$$
\mathbb{E}_\rho(x)=1-\M_n\left(\frac{n}{\omega_n}\int
\frac{\rho(\zeta)d\zeta}{|x-\zeta|^n}\right).
$$
where $\rho$ is a density function.

\begin{cor}
\label{c:bjorn} Function $E_\rho(x)$ is analytic if and only
$\mathbb{E}_\rho(x)$ is. Moreover, these functions linked by the
following identity
\begin{equation}\label{E-E}
\mathbb{E}_\rho(x)=\phi_{2/n}\circ E_\rho(x).
\end{equation}
\end{cor}

\begin{proof}
The cases $n=1$ and $n=2$ are trivial. For $n\geq 3$ we notice
that the desired property follows from (\ref{E-E}) and the fact
that $\phi'_\alpha(0)\ne 0$ (see (\ref{phi-prime}) below).
\end{proof}

\begin{proof}[Proof of Theorem~\ref{theo:phi}]
First we consider \textbf{(i)}. The case $\alpha=1$ is trivial.
Let $\alpha>0$, $\alpha\ne1$ and $F_\alpha(x)$ be the solution to
(\ref{equ:extM}). We notice that this function is determined
uniquely by virtue of the condition $F_\alpha(0)=0$, and it is a
real analytic function of $x$ in $(0,+\infty)$. It follows that
$\phi_\alpha(t)$ also is a real analytic function of $t$ for
$t\in(0,1)$ it is bounded there: $|\phi_\alpha(t)|<1$. Moreover,
$y=\phi_\alpha(t)$ satisfies the following differential equation
\begin{equation}\label{e:diff-phi}
y'(t)=\frac{1-(1-y(t))^\alpha}{\alpha t}, \quad t\in(0,1),
\end{equation}
and the initial condition has to be transformed to $
\phi_\alpha(1)=0.
$

Now we prove that $\phi_\alpha(x)$ admits an analytic continuation
in a small disk in the complex plane. Let us define the following
auxiliary function
$$
S(\z):=\exp\biggl(-\int\limits_{\z}^1 \frac{\alpha d\xi
}{1-(1-\xi)^{\alpha}}\biggr).
$$
Here we fix the branch of $(1-\xi)^{\alpha}$ which assumes the
value $1$ at $\xi=0$. Then $S(\z)$ is a single-valued holomorphic
function  in the unit disk $\DD(1)$, where
$$
\DD(r)=\{\xi\in\Com{}:|\z|<r\}.
$$
Moreover, we have  $S(\z)\ne 0$ and
\begin{equation}\label{e:S}
S'(\z)=\frac{\alpha S(\z)}{1-(1-\z)^\alpha}.
\end{equation}

On the other hand, we notice that the following renormalization of
the above integrand
$$
\frac{\alpha }{1-(1-\xi)^{\alpha}}-\frac{1}{\xi}=\frac{\alpha
\xi-1+(1-\xi)^{\alpha}}{(1-(1-\xi)^{\alpha})\xi}=\frac{\alpha-1}{2}+\frac{(\alpha-1)(2\alpha-1)}{12}\xi+\ldots.
$$
is  an analytic  function in $\DD$ since it admits a regular
Taylor expansion near $\xi=0$. This allows us to rewrite the above
definition of $S(\z)$ as follows
$$
S(\z)=\z\,\exp\biggl\{-\int\limits_{\z}^1 \biggl(\frac{\alpha
}{1-(1-\xi)^{\alpha}}-\frac{1}{\xi}\biggr)d\xi\biggr\}.
$$
In particular, this implies
\begin{equation}\label{e:S1}
c_\alpha:=S'(0)=\exp\biggl(-\int\limits_{0}^1
\biggl(\frac{\alpha}{1-(1-x)^{\alpha}}-\frac{1}{x}\biggr)dx\biggr)\ne
0.
\end{equation}

On the other hand,
$$
\re \frac{S'(\z)\z}{S(\z)}=\re \frac{\z}{1-(1-\z)^{\alpha}}=\re
\frac{1}{\alpha}(1+\frac{\alpha-1}{2}\z+\ldots).
$$
Hence for $r>0$ sufficiently small
\begin{equation}\label{star}
\re \frac{S'(\z)\z}{S(\z)}>0, \quad \z\in\DD(r)
\end{equation}
 Taking into account (\ref{e:S1}),
(\ref{star}), and the well-known Alexander's property
\cite[p.~41]{Dur} we conclude that the function $S(z)$ is starlike
in $\DD(r)$, and therefore univalent there.

Let $\psi(z)$ be the inverse function to $S(z)$. Clearly, it is
defined in some small disk $\DD(\epsilon)$ which is contained in
the image $S(\DD(r))$. Moreover, by its definition $\psi(z)$
assumes real values for real $z\in\DD(\epsilon)$. We also have
\begin{equation}\label{phi-prime}
\psi(0)=0, \qquad \psi'(0)=\frac{1}{S'(0)}=\frac{1}{c_\alpha}.
\end{equation}
Furthermore, differentiation of the identity $S(\psi(z))=z$
together with (\ref{e:S}) yields
$$
1=S'(\psi(z))\psi'(z)=\frac{z\psi'(z)}{1-(1-\psi(z))^\alpha},
$$
consequently,  $y=\psi(z)$ is a solution of (\ref{e:diff-phi}) in
$\DD(\epsilon)$.

Our next step is to prove that $\psi(z)$ is the desired analytic
continuation. One suffices to show that $\psi(x)=\phi_\alpha(x)$
in some open subinterval of $(0,1)$, that in turn, is equivalent
to establishing of the following identity
\begin{equation}\label{e:ident}
F_\alpha(x)=1-\psi(e^{-\alpha x})
\end{equation}
for all $x$ in some interval $\Delta\subset (0,+\infty)$.

Taking into account the above remarks, we note that
$g(x):=1-\psi(e^{-\alpha x})$ is a real-valued solution of
(\ref{equ:extM}) in
$\Delta:=(-\frac{1}{\alpha}\ln\epsilon,+\infty)$. On the other
hand, since $g(x)$ satisfies an obvious inequality $g(x)<1$, and
by virtue of the autonomic character of (\ref{equ:extM}), we
conclude that
$$
g(x)=F_\alpha(x+c), \quad x\in \Delta,
$$
for some constant $c\in\R{}$. Thus, we have only to check that
$c=0$.

To this aim, we note that
$$
\alpha x=\int\limits_{0}^{F_\alpha(x)}\frac{\alpha
dt}{1-t^\alpha}=
-\int\limits_{0}^{F^\alpha_\alpha(x)}\frac{1-\tau^\frac{1-\alpha}{\alpha}}{1-\tau
\phantom{mn}}d\tau +
\int\limits_{0}^{F^\alpha_\alpha(x)}\frac{d\tau}{1-\tau}d\tau.
$$
Since $\lim_{\tau\to\infty}F_\alpha(\tau)=1$, we arrive at
$$
\lim_{x\to+\infty}(1-F^\alpha_\alpha(x))e^{x\alpha}=
\exp\biggl(-\int\limits_{0}^1
\frac{1-\tau^{\frac{1-\alpha}{\alpha}}}{1-\tau\phantom{mn}}\;d\tau\biggr)=
{\alpha\gamma_\alpha},
$$
or
$$
\lim_{x\to+\infty}(1-F_\alpha(x))e^{x\alpha}=\gamma_\alpha.
$$
As a consequence we have,
\begin{equation}\label{e:exp-exp}\begin{split}
\lim_{x\to+\infty}\psi(e^{-x\alpha})e^{x\alpha}&=\lim_{x\to+\infty}(1-g(x))e^{x\alpha}=\\
&=\lim_{x\to+\infty}[1-F_\alpha(x+c)]e^{x\alpha}=e^{-c\alpha}\gamma_\alpha.
\end{split}
\end{equation}
On the other hand,
\begin{equation}\label{e:exp-exp1}
\lim_{x\to+\infty}\psi(e^{-x\alpha})e^{x\alpha}=
\lim_{t\to+0}\frac{\psi(t)}{t}=\psi'(0)=\frac{1}{c_\alpha}.
\end{equation}
Finally, splitting the integral in the definition of $c_\alpha$
and making the change variables $\tau=(1-t)^\alpha$, we obtain
\begin{equation*}
\begin{split}
\ln \frac{1}{c_\alpha}&=\lim_{s\to+0}\int\limits_{s}^1
\biggl(\frac{\alpha}{1-(1-x)^{\alpha}}-\frac{1}{x}\biggr)dx=
\lim_{s\to+0}\biggl(\int\limits^{(1-s)^\alpha}_0
\frac{\tau^{\frac{1-\alpha}{\alpha}}d\tau}{1-\tau}+\ln s\biggr)=\\
&=\lim_{s\to+0}\biggl(-\int\limits^{(1-s)^\alpha}_0
\frac{1-\tau^{\frac{1-\alpha}{\alpha}}}{1-\tau}d\tau-\ln\frac{1-(1-s)^\alpha}{s}\biggr)
=\\
&=-\ln\alpha-\int\limits^{1}_0
\frac{1-\tau^{\frac{1-\alpha}{\alpha}}}{1-\tau}d\tau=\ln\gamma_\alpha.
\end{split}
\end{equation*}
Thus, combining the latter identity with (\ref{e:exp-exp}), and
(\ref{e:exp-exp1}) we obtain $c=0$, which yields (\ref{e:ident})
and the mentioned analytic continuation property follows.

Now we prove \textbf{(ii)} and \textbf{(iii)}. We note that in
view of (\ref{e:diff-phi})
$$
(\alpha t \phi'_\alpha)'=\alpha
(1-\phi_\alpha)^{\alpha-1}\phi'_\alpha= \alpha
\phi'_\alpha\frac{1-\alpha t \phi_\alpha'}{1-\phi_\alpha},
$$
which implies
\begin{equation}\label{e:ode}
t\phi''_\alpha =\phi_\alpha(t\phi''_\alpha+\phi'_\alpha)-\alpha
t\phi'^2_\alpha.
\end{equation}
Setting
\begin{equation}\label{e:con}
\phi_\alpha(t)=\sum_{k=1}^{\infty} a_kt^k
\end{equation}
for the Taylor series of $\phi_\alpha$ around $t=0$ (we recall
that $\phi_\alpha(0)=0$) we obtain after comparison of the
corresponding  coefficients for all $k\geq 2$
\begin{equation*}\label{e:Taylor-coef0}
\begin{split}
a_k&=\frac{1}{k(k-1)}\sum_{\nu=1}^{k-1}a_\nu a_{k-\nu}
[(1+\alpha)\nu-\alpha k]\nu=\sigma_k a_1^k.
\end{split}
\end{equation*}
Here
$$
a_1=\phi'_\alpha(0)=1/c_\alpha=\gamma_\alpha
$$
and $\sigma_k$ are defined as in (\ref{e:Taylor-coef}). This
yields the desired Taylor expansion. Moreover, we show that for
$\alpha\in (0,1)$ the coefficients $\sigma_k>0$ for all $k\geq 1$.
Indeed,
$$
\sigma_k=\frac{1}{2k(k-1)}\sum_{\nu=1}^{k-1}A_{\nu,k-\nu}\sigma_\nu
\sigma_{k-\nu},
$$
where
\begin{equation*}\label{e:Taylor-coef1}
\begin{split}
A_{\nu,k-\nu}=&[(1+\alpha)\nu-\alpha k]\nu+[(1+\alpha)(k-\nu)-\alpha k](k-\nu)\\
=&(1+\alpha)\biggl(\nu-\frac{k}{2}\biggr)^2+\frac{1-\alpha}{2}k^2>0,
\end{split}
\end{equation*}
unless $k=2\nu$ when we also have
$$
A_{\nu,\nu}=2[(1+\alpha)\nu-2\nu\alpha ]\nu=2(1-\alpha)\nu^2>0.
$$

Since $\sigma_1=1$ and for $k\geq 1$ the coefficients
$A_{\nu,k-\nu}$ before $\sigma_\nu \sigma_{k-\nu} $ are positive,
 the positiveness of $\sigma_k$ follows now by induction.

 Thus, $\phi_\alpha(t)$ has the Taylor expansion with positive
 coefficients. By standard facts of the power series theory we
 conclude that the radius $R$ of convergence of (\ref{e:con}) is at
 least $R=1$ since $\phi_\alpha(t)$ is analytic along
 $t\in(-\epsilon,1)$.

 It remains only to prove \textbf{(iv)}. We have $\phi'(0)>0$ which
 yields $\phi_\alpha(t)<0$ for sufficiently small $t<0$. Then a
 standard analysis of (\ref{e:diff-phi}) shows that these property
 holds for \textit{all} negative $t$'s where $\phi_\alpha(t)$ is
 defined. In view of (\ref{e:diff-phi}), this proves the strictly increasing
 character of $\phi_\alpha(t)$.

 In order to prove convexity, we note that (\ref{e:ode}) implies
$$
\phi''_\alpha(t) =\phi'_\alpha(t)\frac{\phi_\alpha(t)-\alpha
t\phi'_\alpha(t)} {t(1-\phi_\alpha(t))}, \quad
\phi''(0)=(1-\alpha)\phi'^2_\alpha(0)>0.
$$
Clearly, it suffices to prove that $\phi''(t)\ne 0$. Assuming the
contradictory, we have $\phi''(t)= 0$ an some point $t\ne0$, and
it follows that
$$
\phi_\alpha(t)-\alpha t\phi'_\alpha(t)=0,
$$
which yields $(1-\phi_\alpha(t))^\alpha=1-\phi_\alpha(t)$. The
contradiction obtained.

Finally, since $\phi_\alpha(t)$ is convex and analytic in its
region of definition, we conclude that it can be infinitely
extended into the left side of $\R{}$.
\end{proof}

\section{Final remark}

 Here we  discuss in short an appearance of the profile function
 $\M_1(w)$
as interpretation of the exponential transform. We recall that the
original result of A.A.~Markov on the $L$-problem asserts that a
sequence of reals $\{s_j\}_{j=0}^{\infty}$ is represented as the
moments
$$
s_k=\int\limits x^k \rho(x)dx
$$
of certain function $0\leq \rho(x)\leq 1$, if and only if there is
a positive measure $d\mu$ such that the following identity holds
\begin{equation*}\label{exp:def}
1-\exp(-\sum_{k=0}^{\infty}\frac{s_k}{z^{k+1}})=\sum_{k=0}^{\infty}\frac{a_k}{z^{k+1}},
\end{equation*}
where
$$
a_k=\int\limits x^kd\mu(x).
$$
For the detailed discussion of this theory see
\cite[p.~72]{Akh-Kr}. The latter moment sequence,
$\{a_j\}_{j=0}^{\infty}$, can be characterized as a standard
\textit{positive} sequence in the sense that the Hankel forms
$$
(a_{i+j})_{i,j=0}^{m}\geq 0
$$
are positive semi-definite for all $m\geq 0$. For simplicity
reasons, we refer to $(s_k)$ as an $L$-sequence.

Given a sequence $(a_{k})_{k=0}^{\infty}$ we set
$$
\widehat{a}(z):=\sum_{k=0}^{\infty}\frac{a_k}{z^{k+1}}
$$
for the corresponding $z$-transform. Our first observation is as
follows.

\begin{prop}
Let $c\in\R{}$, and $\{a_j\}_{j=0}^{\infty}$ and
$\{b_j\}_{j=0}^{\infty}$ be two sequences such that their
generating functions satisfy
\begin{equation}\label{e-ab}
\frac{1}{\widehat{b}(z)}-\frac{1}{\vphantom{\widehat{b}(z)}\widehat{a}(z)}=c.
\end{equation}
Then $\{a_j\}_{j=0}^{\infty}$ is a positive sequence if and only
$\{b_j\}_{j=0}^{\infty}$ is. Moreover, we have
\begin{equation}\label{e:deter}
\det (a_{i+j})_{i,j=0}^{m}= \det (b_{i+j})_{i,j=0}^{m}.
\end{equation}
\end{prop}

\begin{proof}
We prove only (\ref{e:deter}) since it immediately implies the
desired positivity property. Let for definiteness,
$\widehat{a}(z)$ satisfies the positivity condition, i.e. the
corresponding sequence $(a_k)$ is positive semi-definite.

Then the famous result of Stieltjes \cite[Ch.~XI]{Wall} asserts
that given a function $\widehat{a}(z)$ with power series as above,
the following continued $J$-fraction (actually, Jacobi's type)
decomposition holds
\begin{equation}\label{e-aa}
\widehat{a}(z)=\frac{\alpha_0}{\displaystyle\beta_1+z-\frac{\displaystyle\alpha_1}
{\displaystyle\beta_2+z-\frac{\displaystyle
\alpha_2}{\displaystyle\beta_3+z-\ldots}}}.
\end{equation}
Moreover, in this case we have for the determinants
\begin{equation}\label{e:deter1}
\det (a_{i+j})_{i,j=0}^{m}=
\alpha_0^{m+1}\alpha_1^m\alpha_2^{m-1}\cdots
\alpha^2_{m-1}\alpha_m.
\end{equation}
Now, it follows from (\ref{e-ab}) and (\ref{e-aa}) that
\begin{equation*}\label{e-ac}
\widehat{b}(z)=\frac{1}{c +\frac{\displaystyle \vphantom{M^M}
1}{\displaystyle \vphantom{{\widehat{b}(z)}^M}\widehat{a}(z)}}=
\frac{\alpha_0}{\displaystyle
c\alpha_0+\beta_1+z-\frac{\displaystyle\alpha_1}
{\displaystyle\beta_2+z-\frac{\displaystyle
\alpha_2}{\displaystyle\beta_3+z-\ldots}}}.
\end{equation*}
The latter continuous  fraction is the Stieltjes' $J$-fraction for
$\widehat{b}(z)$ and hence we have for its determinants the same
expressions as those in (\ref{e:deter1}), and (\ref{e:deter})
follows.
\end{proof}

\begin{cor}\label{cor:frac}
The sequence $\{s_j\}_{j=0}^{\infty}$ is an $L$-sequence if and
only if
\begin{equation}\label{exp:def1}
\M_1\left(\frac{1}{2}\widehat{s}(z)\right)=\widehat{b}(z)
\end{equation}
for some positive sequence $\{b_j\}_{j=0}^{\infty}$.
\end{cor}

\begin{proof}
Indeed, we have
$$
\M_1(w)=\tanh w= \frac{e^{2w}-1}{e^{2w}+1},
$$
therefore,
$$
\widehat{b}(z)\equiv
\M_1\left(\frac{1}{2}\widehat{s}(z)\right)=\frac{1-v(z)}{1+v(z)},
$$
where $v(z)=\exp(-\widehat{s}(z))$ is the standard exponential
transform of $\widehat{s}(z)$. Since $1-v(z)$ is the generating
function of some positive sequence $(a_k)$, we have
$$
\widehat{b}(z)=\frac{\widehat{a}(z)}{2-\widehat{a}(z)},
$$
or
$$
\frac{1}{\widehat{b}(z)}=\frac{2}{\widehat{a}(z)}-1,
$$
and the required property follows from positivity of
$\widehat{a}(z)/2$.
\end{proof}

\begin{rem}
The previous observation makes it possible to consider an analogue
of the ($n$-dimensional) transform by letting
$$
\mathbb{E}^n_\rho(x):=1-\M_n(V_\rho(x)).
$$
In particular, $ \mathbb{E}^2_\rho(x)= E_\rho(x),
$
while for $n=1$ we have
$$
\mathbb{E}^1_\rho(x)= \frac{2E_\rho(x)}{1+E_\rho(x)}.
$$
\end{rem}


\subsection*{Acknowledgment}
The author is grateful to Bj\"orn Gustafsson, Mihai Putinar and
Serguei Shimorin for conversations crucial to the development of
this paper.
\end{document}